\input amstex
\input psfig.sty
\input Amstex-document.sty

\def\ddt{\frac{d}{dt}}     
\def\nn{\bold n}           
\def\kk{\bold k}           
\def\vv{\bold v}           
\def\RR{\bold R}           
\def\eps{\epsilon}

\def\volumeno{I}
\def\firstpage{525}
\def\lastpage{538}
\pageno=525
\def\shorttitle{Evolution of Curves and Surfaces by Mean Curvature}
\def\shortauthor{B. White}

\topmatter
\title\nofrills{\boldHuge Evolution of Curves and Surfaces by Mean Curvature*}
\endtitle

\author \Large Brian White$^\dag$ \endauthor

\thanks *The preparation of this article was partially funded by NSF grant DMS 0104049.
\endthanks

\thanks $^\dag$Department of Mathematics, Stanford University, Stanford, CA 94305,
USA. \newline E-mail: white\@{}math.stanford.edu, \quad Webpage:
http://math.stanford.edu/\~{}white
\endthanks

\abstract\nofrills \centerline{\boldnormal Abstract}

\vskip 4mm

{\ninepoint This article describes the mean curvature flow, some of the
discoveries that have been made about it, and some unresolved questions.

\vskip 4mm


}
\endabstract
\endtopmatter

\document

\baselineskip 4.5mm \parindent 8mm

\specialhead \noindent \boldLARGE 1. Introduction \endspecialhead

Traditionally, differential geometry has been the study
of curved spaces or shapes
in which, for the most part, time did not play a role.
In the last few decades, on the other hand,
geometers have made great strides in understanding
shapes that evolve in time.
There are many
processes by which a curve or surface can evolve,
but among them one is arguably the most natural:
the mean curvature flow.
This article describes the flow, some of the
discoveries that have been made about it, and some unresolved
questions.

\specialhead \noindent \boldLARGE 2. Curve-shortening flow
\endspecialhead

The simplest case is that of curves in the plane.
Here the flow is usually called the ``curvature flow''
or the ``curve-shortening flow''.
Consider a smooth simple closed curve in the plane,
and let each point move with a velocity equal to the curvature
vector at that point.   What happens to the curve?

The evolution has several basic properties.  First, it makes the curve
smoother.  Consider a portion of a bumpy curve as in figure 1(a).
The portions that stick up move down and the portions that stick down
move up, so the curve becomes smoother or less bumpy as in figure 1(b).
The partial differential
equation for the motion is a parabolic or heat-type equation,
and such smoothing is a general feature of solutions to such equations.
Thus, for example, even if the initial curve is only $C^2$, as it starts
moving it immediately becomes $C^\infty$ and indeed real analytic.

\vskip .2in \vbox{ \centerline{\psfig{figure=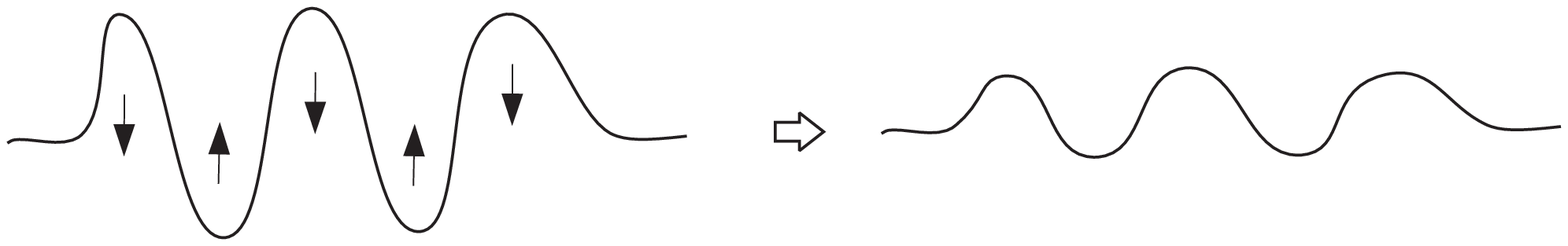,width=4in}}
\centerline{1(a) \hphantom{leave a big chunk of space here} 1(b)}
\centerline{} \centerline{Figure 1: Smoothing} } \vskip 0.1in

However, there is an important caveat:
the smoothing may only be for a short time.
If the curve is $C^2$ at time $0$, it will be real-analytic for times
$t$ in some interval $(0,\eps)$.  But because the equation of motion
is nonlinear, the general theory of parabolic equations does not preclude
later singularities.  And indeed, as we shall see, any curve must eventually
become singular under the curvature flow.

The simplest closed curve is of course a circle.  The flow clearly
preserves the symmetry, so in this case it is easy to solve the equation of
motion explicitly.   One finds that a circle of radius of radius $r$ at time
$0$ shrinks to a circle of radius $\sqrt{r^2-2t}$ at time $t$,
so that time $t=r^2/2$
the circle has collapsed to a point and thereby become singular.

The second fundamental property is that arclength decreases.
The proof is as follows.  For any evolution of curves,
$$
  \ddt \text{(length)} = - \int \kk\cdot \vv \,ds,
$$
where $\kk$ is the curvature vector, $\vv$ is the velocity,
and $ds$ is arclength.  For the curvature flow, $\vv=\kk$, so the right
hand side of the equation is clearly negative.
Indeed, the proof shows that this flow is, in a sense, the gradient
flow for the arclength functional.
Thus, roughly speaking, the curve evolves so as to reduce its
arclength as rapidly as possible.  This explains the name
``curve-shortening flow'', though many other flows also reduce
arc-length.

The third property is that the flow is collision-free:
two initially disjoint curves must remain disjoint.
The idea of the proof is as follows.  Suppose that two initially
disjoint curves, one inside the other, eventually collide.
At the first time $T$ of contact,
they must touch tangentially.
At the point of tangency, the curvature of the inner curve is
greater than or equal to the curvature of the outer curve.
Suppose for simplicity that strict inequality holds.  Then (at the point
of tangency) the inner curve is moving inward faster than the outer curve is.
But then at a slighly earlier time $T-\eps$, the curves would have to
cross each other.  But that contradicts the choice of $T$
(as first time of contact), proving that contact can never occur.

This collision avoidance is a special case of the maximum principle
for parabolic differential equations.
The maximum principle also
implies in the same way
that an initially embedded curve must remain embedded.

The fourth fundamental property is that every curve $\Gamma$
has a finite lifespan.
To undertand why, consider a circle $C$ that contains $\Gamma$ in its
interior.  Let $\Gamma$ and $C$ both evolve.  The circle collapses in a finite
time.  Since the curves remain disjoint, $\Gamma$ must disappear before the
circle collapses.

There is another nice way to see that a curve must become singular
in a finite time:

{\bf Theorem 1.} \it
If $A(t)$ is the area enclosed by the curve at time $t$, then
$A'(t)=-2\pi$ until the curve becomes singular.  Thus a singularity
must develop within time $A(0)/2\pi$.
\rm

{\bf Proof.}
For any evolution,
$$
 A'(t) = \int_{\Gamma(t)}\vv\cdot \nn\,ds
$$
where $\Gamma(t)$ is the curve at time $t$,
$\vv$ is the velocity, $\nn$ is the outward unit normal, and $ds$
is arclength. For the curvature flow,
$\vv=\kk$, so
$$
  A'(t) = \int_{\Gamma}\kk\cdot\nn\,ds
$$
which equals $-2\pi$ by the Gauss-Bonnet theorem.\qed

The first deep theorem about curvature flow was proved by Mike
Gage and Richard Hamilton in 1986 \cite{GH}:

{\bf Theorem 2.} \it
Under the curvature flow, a convex curve remains convex and
shrinks to a point.
Furthermore, it becomes asymptotically circular:
if the evolving curve is dilated to keep the enclosed area constant, then
the rescaled curve converges to a circle.
\rm

This theorem is often summarized by stating that
convex curves shrink
to round points.

The proof is too involved to describe here, but I will point out that
the result is not at all obvious.  Consider for example a long thin
ellipse, with the major axis horizontal.   The curvature is greater at the
ends than at the top and bottom, so intuitively it should become rounder.
But the ends are much farther from the center than the top and bottom are,
so it is not clear that they all reach the center at the same time.  Thus
it is not obvious that the curve collapses to a point rather than a segment.

Indeed, there are natural flows that have all the above-mentioned
basic properties of curvature flow but for which the Gage-Hamilton
theorem fails. Consider for example a curve moving in the
direction of the curvature vector but with speed equal to the cube
root of the curvature. Under this flow, any ellipse remains an
ellipse of the same eccentricity and thus does not become
circular. For this flow, Ben Andrews \cite{A1} has proved that any
convex curve shrinks to an elliptical point.  (See also \cite{AST,
SaT}.) For other flows (e.g. if ``cube root'' is replaced by $r$th
root for any $r>3$), a convex curve must shrink to a point but in
a very degenerate way: if the evolving curve is dilated to keep
the enclosed area constant, then length of the rescaled curve
tends to infinity \cite{A3}. More generally, Andrews has studied
the existence and nonexistence of asymptotic shapes for convex
curves for rather general classes of flows \cite{A2, A3}.

Shortly after Gage and Hamilton proved their theorem,
Matt Grayson proved what is still perhaps
the most beautiful theorem in the subject:

{\bf Theorem 3.} \cite{G1}  \it
Under the curvature flow, embedded curves become convex and thus (by
the Gage-Hamilton theorem) eventually shrink to round points.
\rm

Again, the proof is too complicated to describe here, but let me
indicate why the result is very surprising.  Consider the annular region
between a two concentric circles of radii and $r=1$ and $R=2$.  Form a
curve in this annular region by spiraling inward $n$ times, and then
back out $n$ times to make a closed embedded curve.  Figure 2 shows such
a curve with $n=3/2$, but think of $n$ being very large, say $10^{100}$.
Recall that the curve exists for a time at most $R^2/2=2$.  By Grayson's
theorem, the curve manages, amazingly, to unwind itself and become
convex in this limited time. Incidentally, notice that initially, except
for two very small portions, the curve is not even moving fast: its
curvature is no more than that of the inner circle.

\vskip .2in \vbox{ \centerline{\psfig{figure=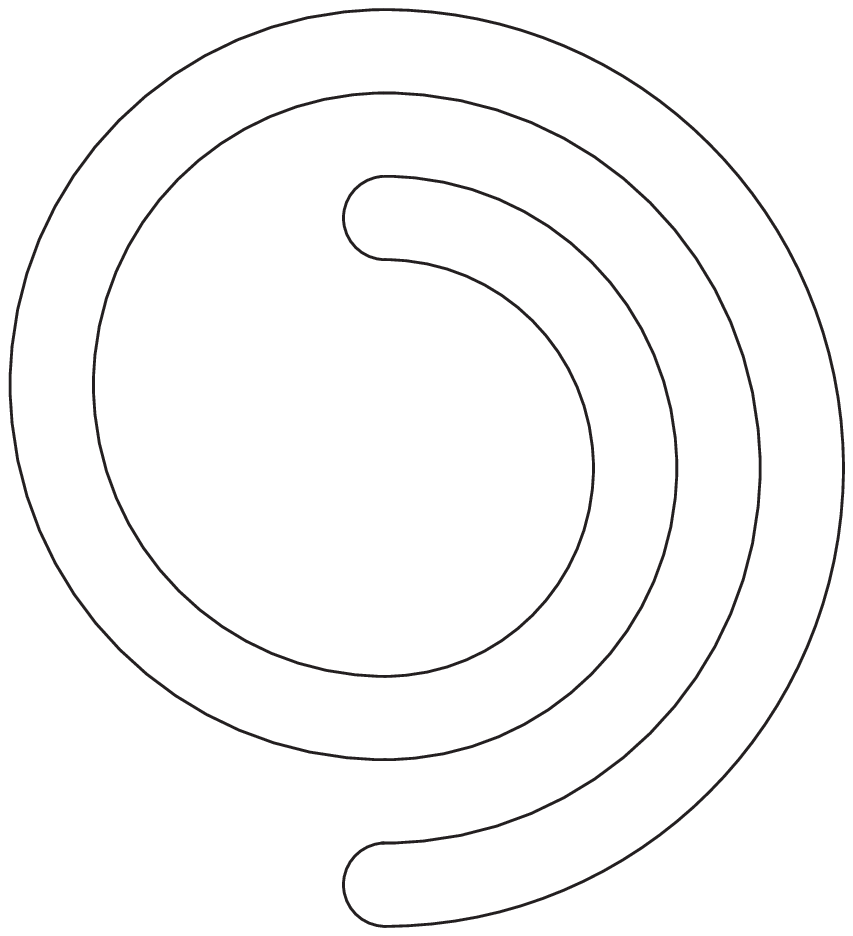,height=1.5in}}
\centerline{} \centerline{Figure 2: Spiral} } \vskip .1in

As a corollary to Grayson's theorem, one gets
an exact formula for the lifespan of any curve.
Recall that the area enclosed by a curve decreases with constant
speed $-2\pi$ as long as the curve is smooth.  By Grayson's theorem,
the curve remains smooth until its area becomes $0$.  Thus the lifespan of any
embedded curve must be exactly equal to the initial area divided by $2\pi$.

Grayson later generalized his theorem by proving that
a closed curve moving on a compact surface by curvature flow must either
collapse to a round point in a finite time or else converge to a simple
closed geodesic as $t\to\infty$ \cite{G2}.

\specialhead \noindent \boldLARGE 3.
       Mean curvature flow for surfaces
\endspecialhead

We now leave curves and consider instead moving surfaces.
Recall that at each point of an $n$-dimensional hypersurface in $\RR^{n+1}$,
there are $n$-principal curvatures $\kappa_1,\kappa_2,\dots,\kappa_n$
(given a choice $\nn$ of unit normal).  Their sum $h$ is
the scalar mean curvature, and the product $H=h\nn$ of the scalar
mean curvature and the unit normal is the mean curvature
vector.  The mean curvature vector does not depend on choice of normal
since replacing $\nn$ by $-\nn$ also changes the sign of the scalar mean
curvature.  In the mean curvature flow, a hypersurface evolves so that
its velocity at each point is equal to the mean curvature vector at that
point.

The basic properties of curvature flow also hold for mean curvature
flow:
\roster
\item Surfaces become smoother for a short time.
\item The area decreases.  Indeed, mean curvature flow may be
  regarded as gradient flow for the area functional.
\item Disjoint surfaces remain disjoint, and embedded surfaces
  remain embedded.
\item Compact surfaces have limited lifespans.
\endroster

The analog of the Gage-Hamilon theorem also holds, as Gerhard Huisken
\cite{H1} proved:

{\bf Theorem 4.} \it
For $n\ge 2$, an $n$-dimensional compact convex surface in $\RR^{n+1}$
must shrink to a round point.
\rm

Oddly enough, Huisken's proof does not apply to the case of curves ($n=1$)
considered by Gage and Hamilton.
Huisken's proof shows that the asymptotic shape is totally umbilic:
at each point $x$, the principal curvatures are all equal
(though {\it a priori} they may vary from point to point).
For $n\ge 2$,
the only totally umbilic surfaces are spheres, but for $n=1$, the condition
is vacuous.

The analog of Grayson's theorem, however, is false for surfaces.
Consider for example two spheres joined by a long thin tube.  The spheres
and the tube both shrink, but the mean curvature along the tube is much
higher than on the spheres, so the middle of the tube collapses down to
a point, forming a singularity.
The surface then separates into two components, which eventually
become convex and collapse to round points.

Thus, unlike a curve, a surface can develop singularities before it
shrinks away.  This raises various questions:
\roster
\item How do singularities affect the subsequent evolution of the surface?
\item How large can the set of singularities be?
\item What is the nature of the singularities?  What does the surface
look like near a singular point?
\endroster

In the rest of this article I will describe some partial answers to these
questions.

A great many results about mean curvature flow have been proved using
only techniques of classical differential geometry and partial differential
equations.  However, most of the proofs are valid only
until the time that singularities first occur.
Once singularities
form, the equation for the flow does not even make sense classically,
so analyzing the flow after a singularity seems to require other techniques.

Fortunately, using non-classical techniques, namely
the geometric measure theory
of varifolds
and/or the theory of viscosity or level-set solutions,
one can define notions of weak solutions for mean
curvature flow and one can prove existence of a solution (with a given initial
surface) up until the time that the surface
disappears.

The definitions are somewhat involved and will not be given here.
The different definitions are equivalent to each other (and to the classical
definition) until singularities form, but are not completely equivalent
in general.
For the purposes of this article, the reader may simply accept that there
is a good way to define mean curvature flow of possibly singular surfaces
and to prove existence theorems.
The notion of mean curvature flow most appropriate to this article
is Ilmanen's ``enhanced Brakke flow of varifolds'' \cite{I}.

\specialhead \noindent \boldLARGE 4.
      Non-uniqueness or fattening
\endspecialhead

If a surface is initially smooth, classical partial differential equations
imply that there is a unique solution of the evolution equation
until singularities form.
However, once a singularity forms, the classical uniqueness theorems
do not apply.
In the early 1990's various researchers, including De Giorgi,
Evans and Spruck, and Chen, Giga, and Goto, asked whether uniqueness could
in fact break down after singularity formation.
They already knew that uniqueness did fail for certain initially singular
surfaces, but they did not know whether an initially regular surface could
later develop singularities that would result in non-uniqueness.

A technical aside: the above-named people did not phrase the question
in terms of uniqueness but rather in terms of ``fattening''.
They were all using a level set or viscosity formulation of mean curvature
flow, in which solutions are almost by definition unique.
But non-uniqueness of the enhanced varifold solutions
corresponds to fattening of the viscosity solution in the following
sense.
If a single initial surface $M$ gives rise to different
enhanced varifold solutions $M^1_t, M^2_t,\dots, M^k_t$, then
the viscosity solution ``surface'' at time $t$ will consist
of the various $M^i_t$'s together with all the points in between. Thus
if $k>1$, the viscosity surface will in fact have an interior.  Since the
surface was initially infinitely thin, in developing an interior it
has thereby ``fattened''.

Recently Tom Ilmanen and I settled this question \cite{IW}:

{\bf Theorem 5.} \it
There is a compact smooth embedded surface in $\RR^3$ for which
uniqueness of (enhanced varifold solutions of) mean curvature evolution
fails.  Equivalently, the viscosity (or level set) solution fattens.
\rm

The idea of the proof is as follows.  Consider a solid torus of revolution
about the z-axis centered at the origin, a ball centered at the origin that
is disjoint from the torus, and $n$ radial segments in the $xy$-plane
joining the ball to the torus.  Call their union $W$.
Now consider a nested one-parameter family of smooth surfaces
$M^\eps$ ($0<\eps<1$) as follows.
When $\eps$ is small, the surface should be a smoothed version
of the set of points at distance $\eps$ from $W$.
This $M^\eps$ looks like a wheel with $n$ spokes.
The portion of the $xy$-plane that is not contained in $M^\eps$
has $n$ simply-connected components, which we regard as holes
between the spokes of the wheel.
As $\eps$ increases, the spokes get thicker and the holes between
the spokes get smaller.
When $\eps$ is close to $1$, the holes should be very small, and
near each hole the surface should resemble a thin vertical tube.

Now let $M^\eps$ flow by mean curvature.  If $\eps$ is small, the
spokes are very thin and will quickly pinch off, separating the
surface into a sphere and a torus.  If $\eps$ is large, the holes
between the spokes are very small and will quickly pinch off, so
the surface becomes (topologically) a sphere. By a continuity
argument, there is an intermediate $\eps$ such that both pinches
occur simultaneously.

For this particular $\eps$, we claim that the simultaneous pinching
immediately results in non-uniqueness, at least if $n$ is sufficiently
large.
Indeed, at the moment of simultaneous pinching, the surface will
resemble a sideways figure 8 curve revolved around the $z$-axis as indicated
in figure 3(a).  Of course the surface will not be fully rotationally
symmetric about the axis, but it will have $n$-fold rotational symmetry,
and here I will be imprecise and proceed as though it were rotationally
symmetric.

\comment \vskip .2in \vbox{
\centerline{\psfig{figure=fatten.eps,width=4in}} \centerline{3(a)
\hphantom{leave a big chunk of space here} ?(b)} \centerline{}
\centerline{Figure ?: Non-uniqueness} } \vskip .5in
\endcomment

\vskip .2in \vbox{ \centerline{\psfig{figure=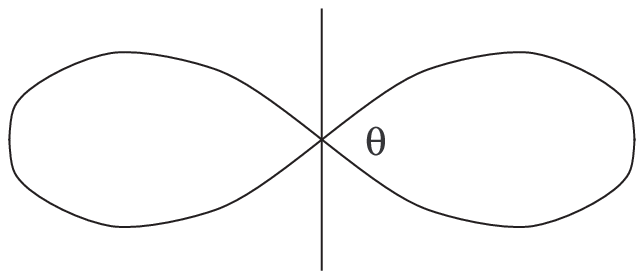,height=1in}}
\centerline{3(a)} \centerline{\psfig{figure=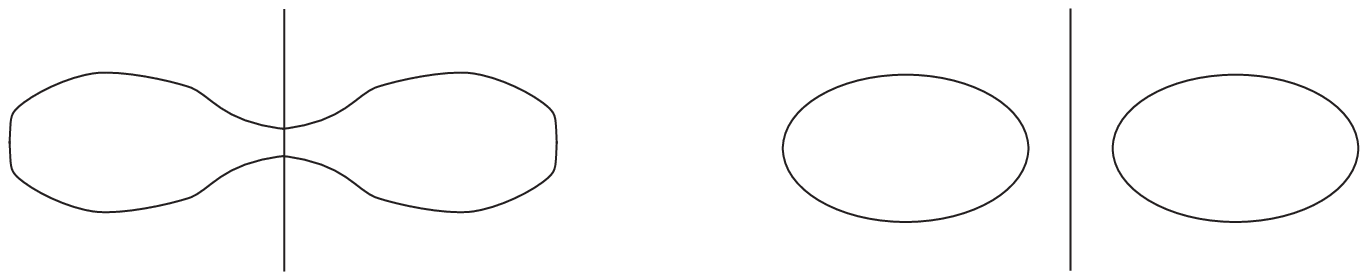,height=1in}}
\centerline{3(b) \hphantom{leave a largish chunk of space here} 3(c)}
\centerline{} \centerline{Figure 3: Non-uniqueness} } \vskip .1in

There is necessarily one evolution in which the surface then becomes
topologically a sphere as in figure 3(b).
If the angle $\theta$ in figure 6 is sufficiently small, there is also
another evolution, in which the surface detaches itself from the $z$-axis
and thereby becomes a torus as in figure 3(c).

One can show that as $n\to\infty$, the angle $\theta$ tends to $0$.
Thus if $n$ is large enough, the angle will be very small and both
evolutions will occur.

The proof unfortunately does not give any bound on how large an $n$
is required.  Numerical evidence \cite{AIT} suggests that $n=4$ suffices.
The case $n=2$ seems to be borderline and the case $n=3$ has not been
investigated numerically.

However, the argument completely breaks down for $n<2$.  Indeed, I would
conjecture that if the initial surface is a smooth embedded sphere
or torus, then uniqueness must hold.

It is desirable to know natural conditions on the initial surface
that guarantee uniqueness.
As just mentioned, genus $\le 1$
may be such a condition.
Mean convexity (described below) and star-shapedness are known to guarantee
uniqueness.  The latter is interesting because the
surface will typically cease to be star-shaped after a finite time,
but its initially starry shape continues to ensure uniqueness \cite{So}.

Fortunately, uniqueness is known to hold generically in a rather strong
sense.  If a family of hypersurfaces foliate an open set in $\RR^{n+1}$,
then uniqueness will hold for all except countably many of the leaves.
Of course any smooth embedded surface is a leaf of such a foliation,
so by perturbing the surface slightly, we get a surface for which uniqueness
holds.

\specialhead \noindent \boldLARGE  5.
   The size of singular sets
\endspecialhead

For general initial surfaces, our knowledge about singular sets is very
limited. Concerning the size of the singular set, Tom Ilmanen \cite{I},
building on earlier deep work of Ken Brakke \cite{BK}, proved the
following theorem.  (See also \cite{ES IV}.)

{\bf Theorem 6.} \it
For almost every initial hypersurface $M_0$ and for almost every time $t$,
the surface $M_t$ is smooth almost everywhere.
\rm

This theorem reminds me of Kurt Friedrichs, who used to say that he
did not like measure theory because when you do measure theory, you have
to say ``almost everywhere'' almost everywhere.

Aside from objections Friedrichs might have had,
the theorem is unsatisfactory in that a much stronger statement should
be true.  But it is a tremendous acheivement and it is the best result to date
for general initial surfaces.

However, for some classes of initial surfaces,
we now have a much better understanding
of singularities.
In particular, this is the case when the initial surface is {\bf mean-convex}.
The rest of this article is about such surfaces.
For simplicity of language,
only two-dimensional surfaces in $\RR^3$ will be discussed, but
the results all have analogs for $n$-dimensional surfaces in $\RR^{n+1}$
or, more generally, in $(n+1)$-dimensional riemannian manifolds.

Consider a compact surface $M$ embedded in $\RR^3$ and bounding
a region $\Omega$.  The surface is said to be ``mean-convex'' if the mean
curvature vector at each point is a nonnegative multiple of the
inward unit normal (that is, the normal that points into $\Omega$.)
This is equivalent to saying that under the mean curvature flow, $M$
immediately moves into $\Omega$.
Mean convexity is a very natural condition for mean curvature flow:
\roster
\item If a surface is initially mean convex, then it remains mean convex
   as it evolves.
\item Uniqueness (or non-fattening) holds for mean convex surfaces.
\endroster
Mean convexity, although a strong condition,
does not preclude interesting singularity formation.
For example, one can connect two spheres by a thin tube as described
earlier in such a way that the resulting surface is mean convex.
Thus neck pinch singularities do occur for some mean convex surfaces.

{\bf Theorem 7.} \cite{W1}   \it
A mean convex surface evolving by mean curvature flow in $\RR^3$
must be completely smooth (i.e., with no singularities) at almost all times,
and at no time can the singular set be more than $1$ dimensional.
\rm

This theorem is in some ways optimal.  For example, consider a torus
of revolution bounding a region $\Omega$.
If the torus is thin enough, it will be mean convex.
Because the symmetry is preserved and because the surface always remains
in $\Omega$, it can only collapse to a circle.  Thus at the time
of collapse, the singular set is one-dimensional.

However, in other ways the result is probably not optimal.
In particular, the result should hold without the
mean convexity hypothesis, and singularities should occur at only
finitely many times.
Indeed, I would conjecture that at each time,
the surface can only have finitely many
singularities unless one or more connected components have collapsed to
curves.  That is, the surface should consist of finitely many connected
components, each of which either is a curve or has only finitely many
singularities.

\specialhead \noindent \boldLARGE  6.
      Nature of mean-convex singularities
\endspecialhead

Recall that when a mean convex surface evolves, it starts moving inward.
Because mean convexity is preserved, it must continue to move inward.
Consequently, the surface at any time lies strictly inside the region
bounded by the surface at any previous time. Since the motion is
continuous and since the surface collapses in a finite time, this
implies that region $\Omega$ bounded by $M_0$ is the disjoint union of
the $M_t$'s for $t>0$. It is convenient and suggestive to speak of the
$M_t$'s forming a foliation of $\Omega$, although it is not quite a
foliation in the usual sense because some of the leaves are singular.
Figure 4 shows the foliation when the initial surface is two spheres
joined by a thin tube.  (The entire foliation is rotationally symmetric
about an axis, so suffices to show the intersection of the foliation
with a plane containing that axis.)

\vskip .2in \vbox{ \centerline{\psfig{figure=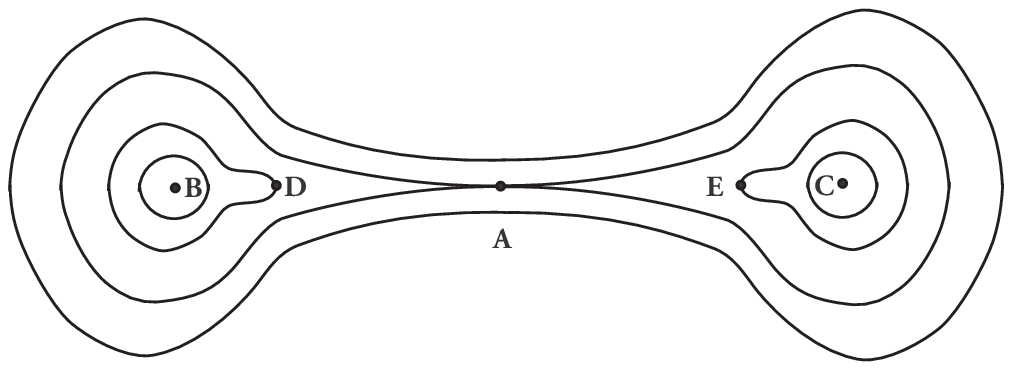,width=4in}}
\centerline{} \centerline{Figure 4: Foliation by evolving mean-convex
surfaces} } \vskip .1in

{\bf Theorem 8.} \cite{W2} \it
Consider a mean convex surface $M_t$ ($t\ge 0$) in $\RR^3$ evolving
by mean curvature flow.  Let $p$ be any singular point in the region $\Omega$
bounded by the initial surface.
If we dilate about $p$ by a factor $\lambda$ and then then let
$\lambda\to \infty$, the dilated foliation must converge subsequentially
to a foliation of $\RR^3$ consisting of
\roster
\item parallel planes, or
\item concentric spheres, or
\item coaxial cylinders.
\endroster
\rm

Let us call such a subsequential limit a {\bf tangent foliation} at $p$.
The tangent foliation consists of parallel planes if and only if $p$
is a regular point (i.e., $p$ has a neighborhood $U$ such that the $M_t\cap U$
smoothly foliate $U$.)
A tangent foliation of concentric spheres corresponds to $M_t$ (or a component
of $M_t$) becoming convex and collapsing as in Huisken's theorem (theorem 4)
to the round point $p$.

In figure 4, the tangent foliation at the ``neck pinch'' point $A$ is
a foliation by coaxial cylinders.  The tangent foliations at points $B$ and
$C$ are by concentric spheres.  All other points are regular points and thus
give rise to tangent foliations consisting of parallel planes.

Incidentally, in cases \therosteritem1 and \therosteritem2, the tangent
foliation is unique.  That is, we have convergence and not just
subsequential convergence in the statement of theorem 8. However, this
is not known in case \therosteritem3.  If one tangent foliation at $p$
consists of cylinders, then so does any other tangent folation at $p$
\cite{S}. But it is conceivable that different sequences of $\lambda$'s
tending to infinity could give rise to cylindrical foliations with
different axes of rotational symmetries. Whether this can actually
happen is a major unsolved problem, exactly analogous to the
long-standing uniqueness of tangent cone problem is minimal surface
theory.

Although the tangent foliation at a singular point carries much information
about the singularity, there are features that it misses.
For example, consider the neck pinch in figure 4, located at the point A.
At a time just after the neck pinch,
the two points $D$ and $E$ on the surface that are nearest to $A$
have very large mean curvature and are therefore moving away from $A$ very
rapidly.
However, such behavior cannot be seen in tangent foliations:
the tangent foliation at any point near $A$ consists of parallel planes,
and the tangent foliation at $A$ consists of coaxial cylinders.

To capture behavior such as the rapid motion away from $A$,
rather than dilating about a fixed point as in
theorem 8, one needs to track a moving point.

{\bf Theorem 9.} \it
Consider a mean convex surface $M_t$ ($t\ge 0$) in $\RR^3$ evolving
by mean curvature flow.
Let $p_i$ be a sequence of points converging to a point $p$ in the region
bounded by $M_0$, and let $\lambda_i$ be a sequence of numbers tending
to infinity.
Translate the $M_t$'s by $-p_i$ and then dilate by $\lambda_i$.
Then the resulting sequence of foliations must converge subsequentially
to a foliation of $\RR^3$ by one of the following:
\roster
\item compact convex sets, or
\item coaxial cylinders, or
\item parallel planes, or
\item non-compact strictly convex surfaces, none of which are singular.
\endroster
The convergence is locally smooth away from the limit foliation's
singular set (a point in case \therosteritem1,
a line in case \therosteritem2, and the empty set
in the other two cases.)
\rm

A foliation obtained in this way is called a {\bf blow-up foliation} at $p$.
Of course if all the $p_i$'s are equal to $p$, we get a tangent foliation
at $p$.

It follows from the smooth convergence that a blow-up foliation is invariant
under mean-curvature flow.  That is, if we let a leaf flow for a time $t$,
the result will still be a leaf.
Consequently, except for the case \therosteritem3 of parallel planes
(which do not move under the flow),
given any two leaves, one will flow to the other in finite time.
Thus we can index the leaves as $M'_t$ in such a way that when $M'_t$ flows
for a time $s$, it becomes $M'_{t+s}$.
In the cases \therosteritem1 and \therosteritem2
of compact leaves and cylindrical leaves,
the indexing interval may be taken to be $(-\infty,0]$.  In case
\therosteritem3,
the indexing interval is $(-\infty,\infty)$.   Thus, except in the case
of parallel planes, the blow-up foliation corresponds to
a semi-eternal or eternal flow of convex sets that sweep out all of $\RR^3$.

As pointed out earlier, blow-up foliations
\therosteritem1, \therosteritem2, and \therosteritem3
already occur as tangent foliations.
(If \therosteritem2 occurs as a tangent foliation, then the compact convex sets
must all be spheres, but in general blow-up foliations, other compact convex
sets might conceivably occur.)

Thus the new case is \therosteritem4.  To see how that case arises, consider
in figure 4 a sequence of points $E_i$ on the axis of rotational symmetry
that converge to the neck pinch point $A$.
Let $h_i$ be the mean curvature at $E_i$ of the leaf of the foliation
that passes through $E_i$.

Now if we translate the foliation by $-E_i$ and then dilate by $\lambda_i=h_i$,
then the dilated foliations converge to a blow-up foliation of $\RR^3$
by convex non-compact surfaces.
Each leaf qualitatively resembles a rotationally
symmetric paraboloid $y=x^2+z^2$.
Furthermore, the leaves are all translates of each
other.   In other words, if we let one of the leaves evolve by mean curvature
flow, then it simply translates with constant speed.

Incidentally, the same points $E_i$
with different choices of $\lambda_i$'s can
give rise to different blow-up foliations.
For if the dilation factors $\lambda_i\to\infty$ quickly compared to $h_i$
(that is, if $\lambda_i/h_i \to \infty$), then the resulting blow-up
foliation consists of parallel planes.
If $\lambda_i\to \infty$ slowly compared to $h_i$
(so that $\lambda_i/h_i\to 0$), then the resulting foliation consists of
coaxial cylinders.

So what is the ``right'' choice of $\lambda_i$?  In a way, it depends
on what one wants to see.   But this example does illustrate
a general principle:

{\bf Theorem 10.} \cite{W2}  \it
Let $p_i\in M_{t(i)}$ be a sequence of regular points converging to a singular
point $p$.  Translate $M_{t(i)}$ by $-p_i$ and the dilate
by $h_i$ (the mean curvature of $M_{t(i)}$ at $p_i$) to get a new
surface $M_i'$.  Then a subsequence of the $M'_i$ will converge
smoothly on bounded subset of $\RR^3$ to a smooth strictly convex
surface $M'$.
\rm

Of course $M'$ is one leaf of the corresponding blow-up foliation.

Theorems 8, 9, and  10 give a rather precise picture of the singular behavior,
but they raise some problems that have not yet been answered:

\roster
\item Classify all the eternal and semi-eternal mean-curvature evolutions
of convex sets that sweep out all of $\RR^3$.
\item Classify those associated with blow-up foliations.
\endroster

The strongest conjecture for \therosteritem2 is that a blow-up foliation
can only consist of planes, spheres, cylinders, or the (unique) rotationally
symmetric translating surfaces.

Many more eternal and semi-eternal evolutions of convex sets are known to
exist.  For instance, given any three positive numbers $a$, $b$, and $c$,
there is a semi-eternal evolution of compact convex sets, each of which
is symmetric about the coordinate planes and one of which cuts off segments
of lengths $a$, $b$, and $c$ from the $x$, $y$, and $z$ axes, respectively.
(This can be proved by a slight modification of the proof given for
example 3 in the ``conclusions'' section of \cite{W2}.)
The case $a=b=c$ is that of concentric spheres, which of course do occur
as a blow-up foliation.
Whether the other cases occur as blow-up foliations is not known.

A very interesting open question in this connection is:
must every eternal evolution of convex sets consist of leaves that move
by translating?  Tom Ilmanen has recently shown that there is a one parameter
family of surfaces that evolve by translation.  At one extreme is the
rotationally symmetric one, which does occur in blow-up foliations.
At the other extreme is the Cartesian product of a certain curve with $\RR$:
$$
 \{ (x,y,z): y = -\ln\cos x,\quad -1<x<1 \}.
$$
This case does not occur in blow-up foliations (\cite{W2}.)

\specialhead \noindent \boldLARGE  7.
         Further reading
\endspecialhead

Three distinct approaches have been very fruitful in
investigating mean curvature: geometric measure theory,
classical PDE, and the theory of level-set or viscosity solutions.
These were pioneered in \cite{BK}; \cite{H1} and \cite{GH};
and \cite{ES} and \cite{CGG} (see also \cite{OS}), respectively.
Surveys emphasizing the classical PDE approach may
be found in \cite{E1} and \cite{H3}.
A very readable and
rather thorough introduction to the classical approach, including
some new results (as well as some discussion of geometric measure theory),
may be found in \cite{E2}.
An introduction to the geometric measure theory
and viscosity approaches is included in \cite{I}.
See \cite{G} for a more extensive introduction to the level set
approach.

Theorems 7, 8, 9, and 10 about mean convex surfaces are from my papers
\cite{W1} and \cite{W2}.
These papers rely strongly on earlier work, for instance
on Brakke's regularity theorem and on
Huisken's monotonicity formula.
Huisken proved theorem 8 much earlier under a hypothesis about the rate
at which the curvature blows up \cite{H2}.
Huisken and Sinestrari \cite{HS 1,2}
independently
proved results very similar to theorems 8, 9, and 10, but only up to the
first occurrence of singularities.

Much of the current interest in curvature flows stems from Hamilton's
spectacular work on the Ricci flow.
For survey articles about Ricci flow, see \cite{CC} and \cite{Ha}.
For discussions of some other interesting geometric flows, see the articles
by Andews \cite{A4} and
by Bray \cite{BH} in these Proceedings.

\specialhead \noindent \boldLARGE
      References
\endspecialhead

\widestnumber\key{CGG}

\ref
\key{A1}
\by B. Andrews
\paper\rm Contraction of convex hypersurfaces by their affine normal
\jour \it J. Differential Geometry,
\vol \rm 43
\yr 1996
\pages 207--230
\endref

\ref
\key{A2}
\bysame
\paper\rm Evolving convex curves
\jour \it Calc. Var. Partial Differential Equations,
\vol \rm 7
\yr 1998
\pages 315--371
\endref

\ref
\key{A3}
\bysame
\paper\rm Non-convergence and instability in the asymptotic behavior
 of curves evolving by curvature
\jour \it Comm. Anal. Geom.,
\vol \rm 10
\yr 2002
\pages 409--449
\endref

\ref
\key{A4}
\bysame
\paper\rm Positively curved surfaces in the three-sphere
\inbook Proceedings of the International Congress of Mathematicians,
   vol. II (Beijing, 2002)
\publ Higher Education Press \publaddr Beijing \yr 2002, 221--230
\endref

\ref
\key{AIC}
\by S. Angenent, T. Ilmanen, and D. L. Chopp
\paper\rm A computed example of nonuniqueness of mean curvature flow in
            $\bold R\sp 3$
\jour \it Comm. Partial Differential Equations,
\vol \rm 20
\yr 1995
\issue 11-12
\pages 1937--1958
\endref

\ref
\key{AST}
\by S. Angenent, G. Sapiro, and A. Tannenbaum
\paper\rm On the affine heat equation for non-convex curves
\jour \it J. Amer. Math. Soc.,
\vol \rm 11
\yr 1998
\issue 3
\pages 601--634
\endref

\ref
\key{BK}
\book The motion of a surface by its mean curvature
\by K. Brakke
\publ Princeton U. Press
\yr 1978
\endref

\ref
\key{BH}
\by H. Bray
\paper\rm Black holes and the penrose inequality in general relativity
\inbook Proceedings of the International Congress of Mathematicians,
  vol. II (Beijing, 2002)
\publ Higher Education Press \publaddr Beijingm \yr 2002, 257--371
\endref

\ref
\key{CC}
\by H.-D. Cao and B. Chow
\paper\rm Recent developments on the Ricci flow
\jour \it Bull. Amer. Math. Soc. (N.S.),
\vol \rm 36
\yr 1999
\issue 1
\pages 59--74
\endref

\ref
\key{CGG}
\paper\rm Uniqueness and existence of viscosity solutions of generalized mean
curvature flow equations
\by Y. G. Chen, Y. Giga, and S. Goto
\jour \it J. Diff. Geom.,
\vol \rm 33
\yr 1991
\pages 749--786
\endref


\ref
\key{E1}
\by K. Ecker
\paper\rm Lectures on geometric evolution equations
\inbook Instructional Workshop on analysis and geometry, part II
    (Canberra, 1995) \yr 1996, 79--107 \publ Austral. Nat. Univ. \publaddr Canberra
\endref

\ref
\key{E2}
\bysame
\paper\rm Lectures on regularity for mean curvature flow
\paperinfo preprint (based on lectures given 2000-2001 at Universit\"at
             Freiburg)
\endref

\ref
\key{ES}
\paper\rm Motion of level sets by mean curvature I
\by L. C. Evans and J. Spruck
\jour \it J. Diff. Geom.,
\vol \rm 33
\yr 1991
\issue 3
\pages 635--681
\moreref
\paper\rm II
\jour \it Trans. Amer. Math. Soc.,
\vol \rm 330
\yr 1992
\issue 1
\pages 321--332
\moreref
\paper\rm III
\jour \it J. Geom. Anal.,
\vol \rm 2
\yr 1992
\issue 2
\pages 121--150
\moreref
\paper\rm IV
\jour \it J. Geom. Anal.,
\vol \rm 5
\issue 1
\yr 1995
\pages 77--114
\endref

\ref
\key{GH}
\paper\rm The heat equation shrinking convex plane curves
\by M. Gage and R. Hamilton
\jour \it J. Differential Geom.,
\vol \rm 23
\yr 1986
\pages 417--491
\endref

\ref
\key{G}
\by Y. Giga
\book Surface evolution equations -- a level set method
\bookinfo Hokkaido U. Tech. Report Series in Math.,
\vol \rm 71
\yr 2002
\publaddr Dept. of Math., Hokkaido Univ., Sapporo 060-0810 Japan
\endref

\ref
\key{Gr1}
\by M. Grayson
\paper\rm The heat equation shrinks embedded plane curves to round points
\jour \it J. Differential Geom.,
\vol \rm 26
\yr 1987
\pages 285--314
\endref

\ref
\key{Gr2}
\bysame
\paper\rm Shortening embedded curves
\jour \it Ann. of Math. (2),
\vol \rm 129
\yr 1989
\issue 1
\pages 71--111
\endref

\ref \key{Ha} \by R. Hamilton \paper\rm    The formation of
singularities in the Ricci flow \inbook Surveys in differential
geometry, Vol.\ II (Cambridge, MA, 1993)  \publ Internat. Press
\publaddr Cambridge, MA \yr 1995, 7--136
\endref

\ref
\key{H1}
\paper\rm Flow by mean curvature of convex surfaces into spheres
\by G. Huisken
\jour \it J. Diff. Geom.,
\vol \rm 20
\yr 1984
\pages 237--266
\endref

\comment
\ref
\key{H2}
\bysame
\paper\rm Asymptotic behavior for singularities of the mean curvature
       flow
\jour \it J. Diff. Geom.,
\vol \rm 31
\yr 1990
\pages 285--299
\endref
\endcomment

\ref
\key{H2}
\paper\rm Local and global behavior of hypersurfaces moving by mean curvature
\bysame
\jour \it Proc. Symp. Pure Math.,
\vol \rm 54
\yr 1993
\publ Amer. Math. Soc.
\pages 175--191
\endref

\ref \key{H3} \bysame \paper\rm Lectures on geometric evolution
equations  \inbook Tsing Hua lectures on geometry and analysis
(Hsinhcu, 1990--1991) \publ Internat. Press \publaddr Cambridge,
MA \yr 1997, 117--143
\endref



\ref
\key{HS1}
\by G. Huisken and C. Sinestrari
\paper\rm Mean curvature flow singularities for mean convex surfaces
\jour \it Calc. Var. Partial Differential Equations,
\vol \rm 8
\yr 1999
\pages 1--14
\endref

\ref
\key{HS2}
\bysame
\paper\rm Convexity estimates for mean curvature flow and singularities
       of mean convex surfaces
\jour \it Acta Math.,
\vol \rm 183
\yr 1999
\pages 45--70
\endref

\ref
\key{I}
\paper\rm Elliptic regularization and partial regularity for
   motion by mean curvature
\by T. Ilmanen
\jour \it Memoirs of the AMS,
\vol \rm 108
\yr 1994
\endref

\ref
\key{IW}
\by T. Ilmanen and B. White
\paper\rm Non-uniqueness for mean curvature flow
 of an initially smooth compact surface
\paperinfo in preparation
\endref

\ref
\key{OS}
\by S. Osher and J. Sethian
\paper\rm Fronts propagating with curvature-dependent speed:
       algorithms based on Hamilton-Jacobi formulations
\jour \it J. Comput. Phys.,
\vol \rm 79
\yr 1988
\pages 12--49
\endref

\ref
\key{SaT}
\by G. Sapiro and A. Tannenbaum
\paper\rm On affine plane curvature evolution
\jour \it J. Funct. Anal.,
\vol \rm 119
\pages 79--120
\issue 1
\endref

\ref
\key{So}
\by H. M. Soner
\paper\rm Motion of a set by the curvature of its boundary
\jour \it J. Differential Equations,
\vol \rm 101
\yr 1993
\issue 2
\pages 313--372
\endref

\ref
\key{St1}
\by A. Stone
\paper\rm A density function and the structure of singularities of
        the mean curvature flow
\jour \it Calc. Var. Partial Differential Equations,
\vol \rm 2
\yr 1994
\pages 443--480
\endref

\comment
\ref
\key{St2}
\bysame
\paper\rm A boundary regularity theorem for mean curvature flow
\jour \it J. Differential Geom.,
\vol \rm 44
\yr 1996
\pages 371--434
\endref
\endcomment

\ref
\key{W1}
\by B. White
\paper\rm The size of the singular set in mean curvature flow
          of mean convex surfaces
\jour \it J. Amer. Math. Soc.,
\vol \rm 13
\yr 2000
\issue 3
\pages 665--695
\endref

\ref
\key{W2}
\bysame
\paper\rm The nature of singularities in mean curvature flow
          of mean convex surfaces
\jour {\it J. Amer. Math. Soc.}
\vol \rm 16
\yr 2003
\pages 123--138
\endref


\enddocument